\documentclass{amsart}

\usepackage[T1]{fontenc}
\usepackage{lmodern}
\usepackage{comment}

\usepackage[francais,english]{babel}

\usepackage{latexsym}
\usepackage{amssymb}

\usepackage[centering,textheight=50.5pc,textwidth=30pc]{geometry}
\usepackage[pdfstartview=XYZ]{hyperref}
\theoremstyle{plain}
\newtheorem{theorem}{Theorem}[section]

\newtheorem{lemma}[theorem]{Lemma}

\theoremstyle{definition}

\theoremstyle{remark}
\newtheorem{remark}[theorem]{Remark}

\theoremstyle{conjecture}

\def\R{\mathbb{R}}

\title[ Ryser's Conjecture under a mild condition]
{A real eigenvector of circulant matrices and a conjecture of Ryser}

\date{\today}

\author[L. H. Gallardo]{Luis H. Gallardo}

\address{University Of Brest,
Mathematics\\
6, Av. Le Gorgeu\\
C.S. 93837\\
29238 Brest Cedex 3, France.}

\email{Luis.Gallardo@univ-brest.fr}



\begin{document}



\begin{abstract}
We prove that there is no circulant Hadamard matrix $H$ with first row $[h_{1},\ldots,h_{n}]$ of order $n>4$, under a condition about a sum of scalar products of
rows of two other circulant matrices of size $n/2$ associated to $H.$ 
\end{abstract}

\maketitle

\section{Introduction}
\label{intro}
A matrix of order $n$ is a  square matrix with $n$ rows.
A \emph{circulant} matrix $A := {\rm{circ}}(a_1,\ldots, a_n)$ of order $n$ is a  matrix of order $n$  of first row $[a_1, \ldots,a_{n}]$ in which each row after the first is
obtained by a cyclic shift of its predecessor by one position. For example, the second row of $A$
is $[a_{n}, a_1,\ldots,a_{n-1}].$
A \emph{Hadamard} matrix $H$ of order $n$ is a matrix of order $n$  with entries in $\{-1,1\}$  such that
$K := \frac{H}{\sqrt{n}}$ is an orthogonal matrix. A \emph{circulant Hadamard} matrix of order $n$ is  a circulant matrix that is Hadamard.  Besides the  two trivial matrices of  order $1$
$H_1 :={\rm{circ}}(1)$ and $ H_2 := -H_1$ the remaining $8$  known circulant Hadamard matrices are $ H_3 :={\rm{circ}}(1,-1,-1,-1), H_4 := -H_3, H_5 :={\rm{circ}}(-1,1,-1,-1), H_6 := -H_5, $ $H_7 := {\rm{circ}}(-1,-1,1,-1), H_8 := -H_7,
 H_9 := {\rm{circ}}(-1,-1,-1,1), H_{10} := -H_9.$

If $H = {\rm{circ}}(h_1,\ldots,h_n)$ is a circulant Hadamard matrix of order $n$ then its \emph{representer} polynomial is the polynomial $R(x):=h_1+h_2x+\cdots+h_nx^{n-1}.$

No one has been able to discover any other circulant Hadamard matrix. Ryser proposed  in $1963$ (see \cite{Ryser}, \cite[p. 97]{Davis}) the conjecture of the non-existence of these matrices
when $n>4.$  Preceding work on the conjecture includes
\cite{EGR,EGRcomb,JedwabLloyd,luisAMEN, luisEJC, Mato, NG,Leung,Turyn}.

The object of the present paper is to prove the conjecture, under a mild condition, in a new special case related to some properties of the real eigenvector $v$ with all entries equal to $1$ of  any circulant matrix. The condition holds for
the $8$ circulant Hadamard matrices of order $4.$ 

In fact the object of the present paper is to prove the following theorem.

\begin{theorem}
\label{HadBM}
Let $H =  \mathrm{circ}(h_1,\ldots,h_n)$ be a circulant Hadamard matrix of order $n \geq 4.$  Let $H_1 := circ(h_1,h_3,h_5,\ldots,h_{n-1})$  and let 
 $H_2 := circ(h_2,h_4,h_6,\ldots,h_{n}).$
Then $n=4$ provided 
\begin{equation}
\label{mainx}
\sum_{j \neq 1, 1 \leq j \leq n/2} \langle\,R_{1},R_{j}\rangle + \sum_{j \neq 1,  1 \leq j \leq n/2} \langle\,S_{1},S_{j}\rangle = \sum_{j \neq 1, 1 \leq j \leq n} \langle\,T_{1},T_{j}\rangle = 0,
\end{equation}
where $R_{j}$, (respectively $S_{j},T_{j}$ is the $j$-th row of the matrix $H_1$ (respectively of the matrices, $H_2$ and $H$) and the  $\langle\,,\rangle$ is the usual scalar product.
\end{theorem}

\begin{remark}
\label{neq4}
When $n=4,$  \eqref{mainx} holds for all $8$ circulant Hadamard matrices $H_3, \ldots,H_{10}$.
\end{remark}

\begin{remark}
\label{apparent}
The condition \eqref{mainx} on Theorem \ref{HadBM} can be proved  {\it heuristically} as follows:
Observe that, by writing explicitly, say, the  first three rows $T_{1},T_{2},$ and $T_{3}$ of $H,$  we  obtain 
$$\langle\,R_{1},R_{2}\rangle + \langle\,S_{1},S_{2}\rangle = \langle\,T_{1},T_{3}\rangle.$$
   But $H$ is Hadamard, thus $\langle\,T_{1},T_{3}\rangle = 0$. Continuing in this manner  we might  eventually obtain the condition (if it is true).
\end{remark}

The necessary tools for the proof of the theorem are given in Section \ref{toolsBM}. The proof of  Theorem  \ref{HadBM} is presented in Section \ref{HaddoneBM},

\section{Tools}
\label{toolsBM}

Our first tool is well -known (\cite{Davis}) and easy to prove.

\begin{lemma}
\label{circulant}
Let $C := circ(c_{1},\ldots,c_{k})$ be a  circulant matrix of order $k$. Let $v := [1, \ldots,1] \in \R^k$. Then,  $v$ is an eigenvector of $C$ 
with asociated eigenvalue $\lambda := c_{1} + \cdots + c_{k}.$
\end{lemma}

The following is well known.  See, e.g., \cite[p. 1193]{HWallis}, \cite[p. 234]{Meisner}, \cite[pp. 329-330] {Turyn}.

\begin{lemma}
\label{regular}
Let $H$ be a regular Hadamard matrix of order $n\geq 4$, i.e., a Hadamard matrix whose row and column sums are all equal. Then $n=4h^2$ for some positive integer $h.$
Moreover, the row and column sums are all equal to $\pm2h$  and each row has $2h^2\pm h$ positive entries and $2h^2 \mp h$ negative entries. Finally, if $H$ is circulant then
$h$ is odd.
\end{lemma}

\begin{lemma}
\label{eigens}
Let $H$ be a circulant Hadamard matrix of order $n,$ let $w= \exp(2 \pi i/n)$ and  let $R(x)$ be its representer polynomial. Then 
\begin{itemize}
\item[(a)]
all the eigenvalues $R(s)$ of $H,$
where $s \in \{1,w,w^2,\ldots,w^{n-1}\},$ satisfy
$$
\vert R(s) \vert = \sqrt{n}.
$$
\item[(b)]
The vector $v := [1, \ldots,1] \in \R^n$
 is eigenvector of $H$ with associated eigenvalue $\lambda = \sqrt{n}.$
\end{itemize}
\end{lemma}

\section{Proof of Theorem \ref{HadBM}}
\label{HaddoneBM}

\begin{proof}

Assume that $n>4.$

By Lemma  \ref{regular} $n$ is even.
Put $$a := \sum_{j \neq 1,  1 \leq j \leq n/2} \langle\,R_{1},R_{j}\rangle \;\text{and} \;  b := \sum_{j \neq 1,  1 \leq j \leq n/2} \langle\,R_{1},R_{j}\rangle.$$  put also 
$\lambda_{1} := h_{1}+h_{3}+ \cdots + h_{n-1}$, the  real eigenvalue of $H_{1}$ associated to the eigevector $v_{0} := [1,\ldots,1] \in \R^{n/2}$, and
$\lambda_{2} := h_{2}+h_{4}+ \cdots + h_{n}$, the  real eigenvalue of $H_{2}$ associated to the same  eigenvector $v_{0}$, (see Lemma \ref{circulant}). Since all $h_j^2=1$ we get 
$ \langle\,R_{1},R_{1}\rangle = n /2 =  \langle\,S_{1},S_{1}\rangle.$ Thus
\begin{equation}
\label{eigen1}
a + n/2 =   \langle\,R_{1}, \sum_{j\, \text{odd}, 1 \leq j \leq n/2} R_{j}\rangle =  \langle\,R_{1}, \lambda_{1} v_{0}\rangle = \lambda_{1}^2.
\end{equation}
and
\begin{equation}
\label{eigen2}
b+ n/2 =   \langle\,S_{1}, \sum_{j \, \text{odd}, 1 \leq j \leq n/2} S_{j}\rangle =  \langle\,S_{1}, \lambda_{2} v_{0}\rangle = \lambda_{2}^2.
\end{equation}

Now, our condition \eqref{mainx} says that

\begin{equation}
\label{coond}
a + b =0
\end{equation}

It follows then  from \eqref{coond} together with  \eqref{eigen1} and \eqref{eigen2} that one has indeed

\begin{equation}
\label{squareseig}
\lambda_1^2 + \lambda_2^2 = n.
\end{equation}

But, it follows from Lemma \ref{eigens} that

\begin{equation}
\label{row1H}
\lambda_{1} + \lambda_{2}=  h_1+ h_2+h_3+ +\cdots + h_n  \in \{\sqrt{n},- \sqrt{n} \}.
\end{equation}

Clearly, we deduce from \eqref{squareseig} and \eqref{row1H} the contradiction

\begin{equation}
\label{clipper}
\lambda_{1}  \lambda_{2}= 0,
\end{equation}
thereby finishing the proof of the theorem.
\end{proof}










\end{document}